\documentclass[12pt]{amsart}
\usepackage{amssymb}
\usepackage[all]{xy}
\usepackage{xcolor}

 \rm

\renewcommand{\(}{\left(}
\renewcommand{\)}{\right)}
\renewcommand{\[}{\left[}
\renewcommand{\]}{\right]}

\renewcommand{\bar}{\overline}
\newcommand{\abs}[1]{\left\lvert#1\right\rvert}

\newcommand{\st}{\:|\:}

\newcommand{\C}{{\mathbb{C}}}
\newcommand{\R}{{\mathbb{R}}}
\newcommand{\Z}{{\mathbb{Z}}}

\renewcommand{\phi}{\varphi}

\renewcommand{\Re}{{\mathrm{Re}}}
\renewcommand{\Im}{{\mathrm{Im}}}

\renewcommand{\H}{{\mathcal{H}}}
\newcommand{\BH}{{\mathcal{B}}(\H)}
\newcommand{\A}{{\mathcal{A}}}
\newcommand{\B}{{\mathcal{B}}}

\theoremstyle{plain}
\newtheorem{thm}{Theorem}
\newtheorem{lem}[thm]{Lemma}

\theoremstyle{definition}

\theoremstyle{remark}

\title{Inductive algebras for the motion group of the plane}

\author{Promod~Sharma}
\author{M.~K.~Vemuri}
\address{Department of Mathematical Sciences\\ IIT (BHU)\\ Varanasi 221 005\\
INDIA}

\begin{document}

\begin{abstract}
Each irreducible representation of the motion group of the plane
has a unique maximal inductive algebra, and it is self adjoint.
\end{abstract}

\keywords{Inductive algebra; Induced representation; Motion group; Lie algebra}
\subjclass[2010]{20C33, 43A80}

\maketitle
\thispagestyle{empty}


\section{Introduction}\label{S:intro}
Let $G$ be a separable locally compact group and $\pi$ an irreducible
unitary representation of $G$ on a separable Hilbert space $\H$.  Let
$\BH$ denote the algebra of bounded operators on $\H$.  An {\em
inductive algebra} is a weakly closed abelian sub-algebra $\A$ of $\BH$
that is normalized by $\pi(G)$, i.e., $\pi(g) \A \pi(g)^{-1} = \A$ for each
$g \in G$.  If we wish to emphasize the
dependence on $\pi$, we will use the term $\pi$-inductive algebra.
A {\em maximal inductive algebra} is a maximal element of the set of
inductive algebras, partially ordered by inclusion.

The identification of inductive algebras can shed light on the
possible realizations of $\H$ as a space of sections of a homogeneous
vector bundle (see e.g. \cite{sharma2022,promod2022,Stegel-2004,Stegel-2006,sl2r,rcr}).
For self-adjoint maximal inductive algebras there is a precise result
known as Mackey's Imprimitivity Theorem, as explained in the
introduction to \cite{Stegel-2004}.  Inductive algebras have also found
applications in operator theory (see e.g. \cite{iahs, Koranyi-2014}).

In Section \ref{S:m2}, we recall the structure and representations of
the motion group of the plane and set up the notation. The main result is proved in the Section \ref{S:m3}.

\section{The Inductive Algebras}\label{S:m2} 
Recall (see \cite[p155]{Sugiura-1975}) that the motion group of the plane is
the semidirect product $$M(2)=SO(2) \ltimes \R^2,$$ with the standard
action of $SO(2)$ on $\R^2$. The group $M(2)$ may be identified with the group
\begin{equation*}
\left\{
\left.
\left[
\begin{matrix}
a & b \\
0 & 1
\end{matrix}
\right] 
\; \right | \:
a, b \in \C,\, \abs{a}=1
\right\}.
\end{equation*}
Let $G=M(2)$. For each $n\in \Z$, let
\begin{equation*}
  \rho_n
  \left[
\begin{matrix}
a & b \\
0 & 1
\end{matrix}
\right] =a^n,
\end{equation*}
then each $\rho_n$ is a character of $G$. These, of course, have the trivial inductive algebra $\C$, which is self adjoint. Let $\mathbb{T}=\{z \in \C \st |z|=1 \}$. Recall (see \cite[p157]{Sugiura-1975}), the irreducible representations other than the characters  of $G$ are parameterized by a non zero real number
$\xi$ and act on the Hilbert space $L^2(\mathbb{T})$ by the formula:
\begin{equation*}
\begin{aligned}
(\pi_\xi\(\[\begin{matrix}a & b \\ 0 & 1  \end{matrix}\]\)F)(z)
=&\; e^{i\xi \Re(b\bar{z})}F(a^{-1}z).
\end{aligned}
\end{equation*}
For each $\phi \in L^{\infty}(\mathbb{T}) $, let $m_\phi: L^2(\mathbb{T})\to L^2(\mathbb{T})$ be
defined by $m_\phi(F)=\phi F$.  Let
\begin{equation*}
\B=\{m_\phi \st \phi \in L^2(\mathbb{T})\}.
\end{equation*}
Then $\B$ is a maximal-abelian subalgebra of 
$\mathcal{B}(L^2(\mathbb{T}))$, and $\B$ is $\pi_{\xi}$-inductive.  Therefore $\B$ is a maximal
$\pi_{\xi}$-inductive algebra.  Moreover, it is self-adjoint.
Our main result is the following theorem.
\begin{thm}
$\B$ is the only maximal $\pi_{\xi}$-inductive algebra.
 \end{thm} 
\section{The proof}\label{S:m3}
Let $\mathfrak{g}$ be the Lie algebra of $G$. Let $\mathcal{H}^\infty$ be the space of all smooth vectors in $\mathcal{H}$ for the representation $\pi_\xi$. Recall that $v\in \mathcal{H}$ is smooth if the map $g \mapsto \pi_\xi(g)v$ is smooth mapping on $G$. The Lie algebra $\mathfrak{g}$ acts on $\mathcal{H}^\infty$ \cite{knapp}.
Observe that $\mathfrak{g}$ is span over $\R$ of the matrices:
\begin{equation*}
L= \[
\begin{matrix}
i & 0\\
0 & 0 
\end{matrix}
\],\quad
M_x= \[
\begin{matrix}
0 & 1\\
0 & 0
\end{matrix}\],\quad \text{and}\quad
M_y= \[
\begin{matrix}
0 & i\\
0 & 0
\end{matrix}
\].
\end{equation*}
Let $e_n(z)=z^n$.
Then
\begin{equation*}
\begin{aligned}
  (d\pi_\xi(L)e_n)(z)  =&\;  \left.\frac{d}{dt}\right|_{t=0}e_n(e^{-it}z)\\
  =&\; -ine_n(z),
\end{aligned}
\end{equation*}

\begin{equation*}
\begin{aligned}
  (d\pi_\xi(M_x)e_n)(z)
=&\;  \left.\frac{d}{dt}\right|_{t=0}e^{i \xi \Re(t\bar{z})}e_n(z)\\
=&\;  i \xi \Re(\bar{z})e_n(z)\\
=&\; i \xi \Re(z)e_n(z).\\
=&\; \frac{i\xi}{2}(e_{n+1}(z)+e_{n-1}(z)), \qquad \text{and}
\end{aligned}
\end{equation*}

\begin{equation*}
\begin{aligned}
(d\pi_\xi(M_y)e_n)(z)
=&\;  \left.\frac{d}{dt}\right|_{t=0}(e^{i \xi \Re(i t\bar{z})}e_n)(z)\\
=&\; i \xi \Re(i\bar{z}) e_n(z)\\
=&\;\xi \Im(z)e_n(z)\\
=&\;\frac{\xi}{2} (e_{n+1}(z)-e_{n-1}(z)).
\end{aligned}
\end{equation*}
Extend $d \pi$ to a $\C$-linear map on $\mathfrak{g} \otimes \C$.
Let $M= M_{y}+i M_{x}$ and $\bar{M}=M_{y}-i M_{x}$.
Then
\begin{equation*}
\begin{aligned}
d\pi_\xi(M)e_n
=&\;  d\pi_\xi(M_y)e_n+i d\pi_\xi(M_x)e_n\\
=&\; - \xi e_{n-1}, \qquad\text{and}\\
d\pi_\xi(\bar{M})e_n
=&\;  \xi e_{n+1}.
\end{aligned}
\end{equation*}

Let
\begin{equation*}
K=\left\{\left.
\begin{bmatrix} w & 0\\ 0 & 1 \end{bmatrix} \:\right|\:
\abs{w}=1
\right\}.
\end{equation*}

Then $K$ is a compact subgroup of $G$.  Let $dk$ denote the Haar
probability measure on $K$.  For $m \in \Z$, define
$\chi_m:K\to\C$ by
\begin{equation*}
\chi_m\(\begin{bmatrix} w & 0\\ 0 & 1 \end{bmatrix}\)=w^m.
\end{equation*}
Let $\mathcal{A}\subseteq \mathcal{B}(\mathcal{H})$ be a $\pi_\xi$-inductive algebra. For each $g\in G$ define a map $\kappa(g): \mathcal{B}(\mathcal{H})\to \mathcal{B}(\mathcal{H})$ by $T \mapsto \pi_\xi(g)T\pi_\xi(g)^{-1}$. Let $\mathcal{A}_m=\{T\in \mathcal{A}~ | ~\kappa(g)T=\chi_m(g)T, g\in K \}$. For $T \in \mathcal{B}(\mathcal{H})$, define
$\Pi_m(T)=\int_K \pi_\xi(k) T \pi_\xi(k)^{-1} \bar{\chi_m(k)} \, dk$, then $\mathcal{A}_m = \Pi_m(\mathcal{A})$. By the Peter-Weyl theorem $\bigoplus_{m\in \Z} \mathcal{A}_m$ is sequentially dense in $\mathcal{A}$ (see \cite{sl2r}).  Note that if $T\in \mathcal{A}_m$ then $Te_n=c_{m+n}e_{m+n}$ for some $c_{m+n} \in \C$. Let $\mathcal{A}^\infty$ be the collection of all smooth vectors in $\mathcal{A}$. We use the following facts from Section 3 of \cite{sl2r}:\\
\begin{itemize}
\item $\mathcal{A}_m\cap \mathcal{A^\infty}$ is sequentially dense in $\mathcal{A}_m$.\\
\item $T\in \mathcal{A}^\infty $ and $ F \in \mathcal{H}^\infty $ then $TF \in \mathcal{H}^\infty$.\\  
\item $\mathcal{A}^\infty$ is invariant under $d\kappa(\mathfrak{g})$ as well as $\kappa(G)$.
\end{itemize}
For $X\in\mathfrak{g}$, let $T_X=[T,d\pi_\xi(X)]$.  If $T\in\mathcal{A}_m \cap \mathcal{A}^{\infty}$, then
\begin{equation*}
\begin{aligned}
T_L(e_n)  =&\;  [T, d\pi_\xi(L)](e_n)\\
=&\;  T(d\pi_\xi(L))(e_n) - d\pi_\xi(L)(Te_n)\\
=&\;  T(-ine_{n})-c_{m+n}d\pi_\xi(L)(e_{n+m})\\
=&\;  -i n c_{m+n}e_{m+n}-c_{m+n }i (m+n) e_{m+n}\\
=&\;  -i n c_n e_{m+n},\\
T_M(e_n)
=&\;  [T, d\pi_\xi(M)](e_n)\\
=&\;   \xi(c_{m+n}-c_{m+n-1})e_{m+n-1},\qquad\text{and}\\
T_{\bar{M}}(e_n)
=&\;  [T, d\pi_\xi(\bar{M})](e_n)\\
=&\; \xi(c_{n+m+1}-c_{n+m})e_{n+m+1}.
\end{aligned} 
\end{equation*}
\begin{lem}\label{L:A0}
$\mathcal{A}_0 = \C I$, where $I$ is the identity operator.
\end{lem}
\begin{proof}
It is clear that $\C I \subset \mathcal{A}_0 $. Let $T \in \mathcal{A}_0$. Then
there exists $c_n\in\C$ such that $Te_n=c_ne_n$  for $n\in \Z$. Therefore 
\begin{equation*}
\begin{aligned}
0
=&\; [T,T_M](e_n) \\
=&\; i\xi(c_{n+1}-c_{n})^2e_{n+1}.
\end{aligned}
\end{equation*}
It follows that  $c_{n+1}=c_{n}$ for all $n \in \mathbb{Z}$.
\end{proof}

Define $S:\mathcal{H}\to \mathcal{H}$ by the formula
\begin{equation*}
S(e_n)=e_{n+1}, \qquad  n \in \mathbb{Z}.
\end{equation*}
\begin{lem}
$\mathcal{A}_m = \C S^m$ for all $m\in \Z$.
\end{lem}

\begin{proof}
The case $m=0$ is Lemma \ref{L:A0}.
Let $T \in \mathcal{A}_1$. Then there exists $c_n\in \C$ such that $Te_n=c_{n+1}e_{n+1}$  for $n\in\Z$ . Since  $T_M \in \mathcal{A}_{0}$, there exists $k\in \C$
 such that $T_M=kI$. Therefore $T_Me_n=ke_n$
for all $n \in \Z$.  Therefore
\begin{equation*}
\begin{aligned}
\xi(c_{1+n}-c_n)e_n = & \; ke_n ,\\
\therefore \xi(c_{1+n}-c_n) = & \; k  ,\\
\therefore c_{1+n} = & \;c_n+\frac{k}{\xi},  \qquad\forall n \in \Z. \\
\end{aligned}
\end{equation*}

Therefore  $\{c_n\}$ is an arithmetic progression,
i.e.\ there exists $a\in \C$ such that
$c_n=a+\frac{k}{\xi}n$ for all $n \in \Z$.
Since $T$ is a bounded operator, it must be the case that $k=0$.
Thus there exists $a\in \C$ such that $c_n=a$ for all $n\in \Z$. Therefore $T=aS$.
This proves the statement for $m=1$.
Now assume that the statement is true for $m=N$.
Let $T \in \mathcal{A}_{N+1}$. Then there exists $c_n\in \C$  such that $Te_n=c_{N+n+1}e_{N+n+1}$ for $n\in\Z$.
Since $T_M \in \mathcal{A}_N$, there exists $k\in \C$ such that
$T_Me_n=k e_{n+N}$ for all $n\in \Z$. Therefore
\begin{equation*}
\begin{aligned}
\xi(c_{N+n+1}-c_{N+n})e_{N+n} =&\; ke_{N+n},\\
\therefore \xi(c_{N+n+1}-c_{N+n}) =&\; k,\\
\therefore c_{N+n+1} =&\; \frac{k}{\xi}+c_{N+n}, \qquad \forall n\in \Z. 
\end{aligned}
\end{equation*}
Therefore $\{c_n\}$ is an arithmetic progression, i.e.\ there exists
$a\in \C$ such that $c_{n}=a+\frac{k}{\xi}n$ for all $n \in \Z$. Since $T$
is a bounded operator, it must be the case that $k=0$.  Therefore
there exists $a\in \C$ such that $c_{n}=a$ for all $n \in \Z$. Therefore $T=aS^{N+1}$. Therefore $\mathcal{A}_{N+1}=\C S^{N+1}$. By induction it follows that
$\mathcal{A}_m=\C S^m$ for $m>0$. The case $m<0$ is similar.
\end{proof}
\bibliographystyle{amsplain}
\bibliography{v17-iaftmgop.bib}

\end{document}